\documentstyle[11pt,twoside]{article}
\include{graphicx}
\title{A note on a modified Bessel function integral}
\author{\sc R. B.\ Paris \\
{\em Division of Computing and Mathematics}, \\
{\em University of Abertay Dundee, Dundee DD1 1HG, UK}
}
\begin{document}
\def\f#1#2{\mbox{${\textstyle \frac{#1}{#2}}$}}
\def\dfrac#1#2{\displaystyle{\frac{#1}{#2}}}
\def\boldal{\mbox{\boldmath $\alpha$}}
{\newcommand{\Sgoth}{S\;\!\!\!\!\!/}
\newcommand{\bee}{\begin{equation}}
\newcommand{\ee}{\end{equation}}
\newcommand{\lam}{\lambda}
\newcommand{\ka}{\kappa}
\newcommand{\al}{\alpha}
\newcommand{\fr}{\frac{1}{2}}
\newcommand{\fs}{\f{1}{2}}
\newcommand{\g}{\Gamma}
\newcommand{\br}{\biggr}
\newcommand{\bl}{\biggl}
\newcommand{\ra}{\rightarrow}
\newcommand{\mbint}{\frac{1}{2\pi i}\int_{c-\infty i}^{c+\infty i}}
\newcommand{\mbcint}{\frac{1}{2\pi i}\int_C}
\newcommand{\mboint}{\frac{1}{2\pi i}\int_{-\infty i}^{\infty i}}
\newcommand{\gtwid}{\raisebox{-.8ex}{\mbox{$\stackrel{\textstyle >}{\sim}$}}}
\newcommand{\ltwid}{\raisebox{-.8ex}{\mbox{$\stackrel{\textstyle <}{\sim}$}}}
\renewcommand{\topfraction}{0.9}
\renewcommand{\bottomfraction}{0.9}
\renewcommand{\textfraction}{0.05}
\newcommand{\mcol}{\multicolumn}
\date{}
\maketitle
\pagestyle{myheadings}
\markboth{\hfill \sc R. B.\ Paris  \hfill}
{\hfill \sc  A Bessel function integral\hfill}
\begin{abstract}
We investigate the integral 
\[\int_0^\infty \cosh^\mu\!t\,K_\nu(z\cosh t)\,dt \qquad \Re(z)>0,\]
where $K$ denotes the modified Bessel function, for non-negative integer values  of the parameters $\mu$ and $\nu$. When the integers are of different parity, closed-form expressions are obtained in terms of $z^{-1}e^{-z}$ multiplied by a polynomial in $z^{-1}$ of degree dependent on the sign of $\mu-\nu$.
\vspace{0.4cm}

\vspace{0.3cm}

\end{abstract}

\vspace{0.3cm}

\noindent $\,$\hrulefill $\,$

\vspace{0.2cm}

\begin{center}
{\bf 1. \  Introduction}
\end{center}
\setcounter{section}{1}
\setcounter{equation}{0}
\renewcommand{\theequation}{\arabic{section}.\arabic{equation}}
An integral arising in neutron scattering \cite{B,BYR} is given by
\bee\label{e11}
{\cal F}(\mu,\nu;z):=\int_0^\infty \cosh^\mu\!t\,K_\nu(z\cosh t)\,dt \qquad \Re(z)>0,
\ee
where $K_\nu(z)$ is the modified Bessel function of the second kind with the parameters $\mu\geq 0$ and $\nu\in C$. 
The asymptotic expansion of ${\cal F}(\mu, \nu;z)$ for $|z|\ra\infty$ in $\Re (z)>0$ has been obtained by Birrell in \cite{B}.
The integral (\ref{e11}) can also be expressed in terms of a combination of three ${}_1F_2(z^2/4)$ generalised hypergeometric functions. The large-$z$ asymptotics could therefore be derived by means of the asymptotic theory of integral functions of hypergeometric type; see, for example, the summary presented in \cite[Section 2.3]{PK}.

Birrell \cite{B} placed particular emphasis on the special case when $\mu$ and $\nu$ are positive integers of different parity, where ${\cal F}(\mu, \nu;z)$ reduces to a finite expansion in inverse powers of $z$ multiplied by $z^{-1}e^{-z}$.
The approach adopted in \cite{B} was the use of a recurrence relation satisfied by the integral in (\ref{e11}) combined with direct evaluation when $\mu=0$, $\nu$ odd and $\mu=1$, $\nu$ even. In this note, we consider the case of integer $\mu$ and $\nu$ from a different viewpoint.
\vspace{0.6cm}

\begin{center}
{\bf 2. \ Evaluation for non-negative integer $\mu$ }
\end{center}
\setcounter{section}{2}
\setcounter{equation}{0}
\renewcommand{\theequation}{\arabic{section}.\arabic{equation}}
We consider the evaluation of ${\cal F}(\mu,\nu;z)$ when $\mu$ and $\nu$ are non-negative integers of different parity by using a well-known property of spherical Bessel functions.
\vspace{0.2cm}

\noindent{2.1\ \ \it The case of even $\mu$}

\noindent
Let us first consider $\mu=2n$ and $\nu=2m+1$, where $m$ and $n$ are non-negative integers. Then we have
\cite[p.~136]{JD}
\[\cosh^{2n}\!t=\frac{1}{2^{2n-1}}\sum_{k=0}^n{}\!' \bl(\!\!\begin{array}{c}2n\\k\end{array}\!\!\br) \cosh 2(n-k)t,\]
where the prime after the summation sign indicates that the term $k=n$ is multiplied by the factor $\fs$. We also make use of the result \cite[p.~253]{DLMF}
\[\int_0^\infty \cosh (a-b)t\,K_{a+b}(2x\cosh t)\,dt=\fs K_a(x) K_b(x) \qquad (\Re (x)>0;\ a, b\in C).\]
It therefore follows that
\bee\label{e20}
{\cal F}(2n, 2m+1;z)=\frac{1}{2^{2n}}\sum_{k=0}^n{}\!' \bl(\!\!\begin{array}{c}2n\\k\end{array}\!\!\br)
K_{m+n-k+\frac{1}{2}}(\fs z) K_{m-n+k+\frac{1}{2}}(\fs z),
\ee
where the Bessel functions are of semi-integer order and so have finite expressions given by \cite[p.~264]{DLMF}
\[K_{s+\frac{1}{2}}(x)=\sqrt{\frac{\pi}{2x}}\,e^{-x} \sum_{k=0}^s\frac{(s+k)!}{k! (s-k)!}\,(2x)^{-k}\qquad (s=0, 1, 2, \ldots).\]

Then, when $m\geq n$ we obtain
\begin{eqnarray}
{\cal F}(2n, 2m+1;z)&=&\frac{\pi e^{-z}}{2^{2n} z}\sum_{k=0}^n{}\!' \bl(\!\!\begin{array}{c}2n\\k\end{array}\!\!\br)\sum_{r=0}^{m+n-k}a_r z^{-r} \sum_{s=0}^{m-n+k}b_s z^{-s}\label{e22}\\
&=&\sum_{k=0}^n{}\!' \bl(\!\!\begin{array}{c}2n\\k\end{array}\!\!\br) \sum_{p=0}^{2m} c_p(k) z^{-p},\nonumber
\end{eqnarray}
where
\bee\label{e22a}
a_r=\frac{(m+n-k+r)!}{r! (m+n-k-r)!},\qquad b_s=\frac{(m-n+k+s)!}{s! (m-n+k-s)!}
\ee
and we have defined 
\begin{eqnarray*}
c_p(k)&\equiv& c_p(k;m,n)=\sum_{r+s=p} a_r b_s\\
&=&\sum_{r=0}^p\frac{(m\!+\!n\!-\!k\!+\!r)! (m\!-\!n\!+\!k\!+\!p\!-\!r)!}{r! (p\!-\!r)! (m\!+\!n\!-\!k\!-\!r)! (m\!-\!n\!+\!k\!-\!p\!+\!r)!}.
\end{eqnarray*}

We therefore find the result
\bee\label{e21}
{\cal F}(2n, 2m+1;z)=\frac{\pi e^{-z}}{2^{2n} z} \sum_{p=0}^{2m} C_p(m,n) z^{-p} \qquad (m\geq n),
\ee
where the coefficients $C_p(m,n)$ are given by
\bee\label{e21a}
C_p(m,n)=\sum_{k=0}^n{}\!' \bl(\!\!\begin{array}{c}2n\\k\end{array}\!\!\br) c_p(k)\qquad (0\leq p\leq 2m).
\ee

In the case $m<n$, we can proceed in a similar manner but now care has to be taken with the second Bessel function appearing in (\ref{e20}). Making use of the fact that $K_{-\nu}(x)=K_\nu(x)$, we write
\bee\label{e25a}
K_{m-n+k+\frac{1}{2}}(\fs z)=\left\{\begin{array}{ll} K_{n-m-k-\frac{1}{2}}(\fs z) & \mbox{for}\ \ k\leq n\!-\!m\!-\!1\\
\\
K_{m-n+k+\frac{1}{2}}(\fs z) & \mbox{for}\ \ k\geq n\!-\!m. \end{array}\right.
\ee
Then the product in (\ref{e22}) becomes
\[\sum_{r=0}^{m+n-k}a_r z^{-r} \sum_{s=0}^S b_s z^{-s},\qquad S=\left\{\begin{array}{ll}n\!-\!m\!-\!k\!-\!1 &\mbox{for}\ \ k\leq n\!-\!m\!-\!1\\m\!-\!n\!+\!k&\mbox{for}\ \ k\geq n\!-\!m\end{array}\right.\]
to yield
\[{\cal F}(2n,2m+1;z)=\frac{\pi e^{-z}}{2^{2n} z}\sum_{k=0}^n{}\!' \bl(\!\!\begin{array}{c}2n\\k\end{array}\!\!\br)
\sum_{p=0}^{2n-1} d_p(k) z^{-p}.\]
The coefficients $d_p(k)\equiv d_p(k;m,n)$ are defined by
\bee\label{e24b}
d_p(k)=\sum_{r+s=p}a_r {\hat b}_s,\qquad {\hat b}_s=\left\{\begin{array}{ll}
\dfrac{(n\!-\!m\!-\!k\!-\!1\!+\!s)!}{s! (n\!-\!m\!-\!k\!-\!1\!-\!s)!} & \mbox{for}\ \ k\leq n\!-\!m\!-\!1\\
\\
\dfrac{(m\!-\!n\!+\!k\!+\!s)!}{s! (m\!-\!n\!+\!k\!-\!s)!} & \mbox{for}\ \ k\geq n\!-\!m,\end{array}\right.
\ee
where the $a_r$ are as specified in (\ref{e22a}). 

Consequently, we obtain the result
\bee\label{e23}
{\cal F}(2n, 2m+1; z)=\frac{\pi e^{-z}}{2^{2n} z}\sum_{p=0}^{2n-1} D_p(m,n) z^{-p}\qquad (m<n),
\ee
where the coefficients $D_p(m,n)$ are given by
\bee\label{e23b}
D_p(m,n)=\sum_{k=0}^n{}\!' \bl(\!\!\begin{array}{c}2n\\k\end{array}\!\!\br) d_p(k)\qquad (0\leq p\leq 2n-1).
\ee
\vspace{0.2cm}

\noindent{\it 2.2\ \ The case of odd $\mu$}

\noindent
When $\mu=2n+1$, we can make use of the expansion \cite[p.~136]{JD}
\[\cosh^{2n+1}\!t=\frac{1}{2^{2n}}\sum_{k=0}^n \bl(\!\!\begin{array}{c}2n+1\\k\end{array}\!\!\br)\cosh (2n-2k+1)t\]
to obtain
\bee\label{e25}
{\cal F}(2n+1,2m;z)=\frac{1}{2^{2n+1}} \sum_{k=0}^n \bl(\!\!\begin{array}{c}2n+1\\k\end{array}\!\!\br)
K_{m+n-k+\frac{1}{2}}(\fs z) K_{m-n+k-\frac{1}{2}}(\fs z).
\ee

Following the same procedure as in the even $\mu$ case, we obtain when $m>n$
\bee\label{e24}
{\cal F}(2n+1,2m;z)=\frac{\pi e^{-z}}{2^{2n+1}z} \sum_{k=0}^n \bl(\!\!\begin{array}{c}2n+1\\k\end{array}\!\!\br)
\sum_{p=0}^{2m-1} {\hat C}_p(m,n) z^{-p}\qquad (m>n),
\ee
where
\bee\label{e24a}
{\hat C}_p(m,n)=\sum_{k=0}^n\bl(\!\!\begin{array}{c}2n+1\\k\end{array}\!\!\br) {\hat c}_p(k)\qquad (0\leq p\leq 2m-1),
\ee
with
\[{\hat c}_p(k)=\sum_{r+s=p}\frac{(m\!+\!n\!-\!k\!+\!r)! (m\!-\!n\!+\!k\!-\!1\!+\!s)!}{r! s! (m\!+\!n\!-\!k\!-\!r)! (m\!-\!n\!+\!k\!-\!1\!-\!s)!}.\]

When $m\leq n$, the second Bessel function in (\ref{e25}) is written as in (\ref{e25a}) with $k$ replaced by $k-1$ to find
\bee\label{e26}
{\cal F}(2n+1,2m;z)=\frac{\pi e^{-z}}{2^{2n+1}z} \sum_{k=0}^{n} \bl(\!\!\begin{array}{c}2n+1\\k\end{array}\!\!\br)
\sum_{p=0}^{2n} {\hat D}_p(m,n) z^{-p}\qquad (m\leq n),
\ee
where
\bee\label{e25a}
{\hat D}_p(m,n)=\sum_{k=0}^n\bl(\!\!\begin{array}{c}2n+1\\k\end{array}\!\!\br) {\hat d}_p(k)\qquad (0\leq p\leq 2n),
\ee
with ${\hat d}_p(k)$ as defined in (\ref{e24b}) with $k\ra k-1$.

We summarise the results obtained in the following theorem.
\newtheorem{theorem}{Theorem}
\begin{theorem}$\!\!\!.$\ \ Let $m$ and $n$ be non-negative integers and $\Re (z)>0$. Then
\bee\label{e27}
{\cal F}(2n,2m+1;z)=\frac{\pi e^{-z}}{2^{2n} z} \,P(z),\qquad {\cal F}(2n+1,2m;z)=\frac{\pi e^{-z}}{2^{2n+1} z} \,{\hat P}(z),
\ee
where $P(z)$ and ${\hat P}(z)$ are polynomials in $z^{-1}$ given by
\[P(z)=\left\{\begin{array}{ll}\displaystyle{\sum_{p=0}^{2m} C_p(m,n)z^{-p}} & (m\geq n)\\
\displaystyle{\sum_{p=0}^{2n-1}D_p(m,n)z^{-p}} &  (m<n)\end{array}\right.,\qquad
{\hat P}(z)=\left\{\begin{array}{ll}\displaystyle{\sum_{p=0}^{2m} {\hat C}_p(m,n)z^{-p}} & (m>n)\\
\displaystyle{\sum_{p=0}^{2n}{\hat D}_p(m,n)z^{-p}} &  (m\leq n).\end{array}\right.\]
The coefficients $C_p(m,n)$, $D_p(m,n)$, ${\hat C}_p(m,n)$ and ${\hat D}_p(m,n)$ are defined in (\ref{e21a}), (\ref{e23b}), (\ref{e24a}) and (\ref{e25a}).
\end{theorem}
\vspace{0.6cm}

\begin{center}
{\bf 3. \  Concluding remarks}
\end{center}
\setcounter{section}{3}
\setcounter{equation}{0}
\renewcommand{\theequation}{\arabic{section}.\arabic{equation}}
As an example, when $\mu=4$, $\nu=7$ ($m=3$, $n=2$) and $\mu=4$, $\nu=3$ ($m=1$, $n=2$) the values of the coefficients $C_p(m,n)$ are
\[8,\ \ 208,\ \ 2520,\ \ 17880,\ \ 76800,\ \ 184320,\ \ 184320\qquad  (m=3, n=2;\ \ 0\leq p\leq 2m),\]
\[8,\ \ 48,\ \ 120,\ \ 120\qquad (m=1, n=2;\ \ 0\leq p\leq 2n-1);\]
whence we find
\[{\cal F}(4,7;z)=\frac{\pi e^{-z}}{2z}\,(1+26z^{-1}+315z^{-2}+2235z^{-3}+9600z^{-4}+23040z^{-5}+23040z^{-6})\]
and 
\[{\cal F}(4,3;z)=\frac{\pi e^{-z}}{2z}\,(1+6z^{-1}+15z^{-2}+15z^{-3}).\]
We note from Theorem 1 that the finite series in $z^{-1}$ terminate at different $p$ index. For even $\mu$, the degree of the polynomial $P(z)$ is $2m$ when $m\geq n$ and $2n-1$ when $m<n$. For odd $\mu$, the degree of the polynomial ${\hat P}(z)$ is $2m$ when $m>n$ and $2n$ when $m\leq n$. This was observed in \cite{B}, although the case $m=n$ was not mentioned.

Since it is easily verified that $c_0(k)=d_0(k)={\hat c}_0(k)={\hat d}_0(k)=1$, the leading coefficients in
the expansions in (\ref{e27}) satisfy
\[C_0(m,n)=D_0(m,n)=\sum_{k=0}^n{}\!' \bl(\!\!\begin{array}{c}2n\\k\end{array}\!\!\br)=2^{2n}-\frac{n!}{2(\fs)_n}\, \bl(\!\!\begin{array}{c}2n\\n\end{array}\!\!\br),\]
\[{\hat C}_0(m,n)={\hat D}_0(m,n)=\sum_{k=0}^n\bl(\!\!\begin{array}{c}2n\!+\!1\\k\end{array}\!\!\br)=2^{2n},\]
which are independent of $m$.

Finally, we mention that the integral
\[{\cal G}(2n,\nu;z)=\int_0^\infty \sinh^{2n}\!t\,K_\nu(z\cosh t)\,dt\]
when $\nu=2m+1$ can be dealt with in the same manner. Since \cite[p.~136]{JD}
\[\sinh^{2n}\!t=\frac{1}{2^{2n-1}} \sum_{k=0}^n{}\!'(-)^k \bl(\!\!\begin{array}{c}2n\\k\end{array}\!\!\br) \cosh 2(n-k)t,\]
it follows immediately from Section 2.1 that
\[{\cal G}(2n,2m+1;z)=\frac{\pi e^{-z}}{2^{2n} z}\,P(z),\]
where the polynomial $P(z)$ is as specified in Theorem 1, with the coefficients $C_p(m,n)$ and $D_p(m,n)$ now defined by
\[C_p(m,n)=\sum_{k=0}^n{}\!'(-)^k \bl(\!\!\begin{array}{c}2n\\k\end{array}\!\!\br) c_p(k),\qquad
D_p(m,n)=\sum_{k=0}^n{}\!'(-)^k \bl(\!\!\begin{array}{c}2n\\k\end{array}\!\!\br) d_p(k).\]
The leading coefficients are again independent of $m$ and given by
\[C_0(m,n)=D_0(m,n)=\sum_{k=0}^n{}\!'(-)^k \bl(\!\!\begin{array}{c}2n\\k\end{array}\!\!\br)=0;\]
this means that a $z^{-1}$ can be factored out leaving the polynomial $P(z)$ with degree reduced by unity.

\end{document}